\documentclass[12pt,reqno]{amsart}
\usepackage[english]{babel}
\usepackage[utf8]{inputenc}

\usepackage{amsmath, amsthm, amscd, amsfonts, amssymb, graphicx, color, mathrsfs}
\usepackage[bookmarksnumbered, colorlinks, plainpages]{hyperref}
\usepackage[all]{xy}
\usepackage[utf8]{inputenc}
\usepackage[english]{babel}
\usepackage{slashed}
\usepackage{cleveref}
\usepackage{mathtools}
\usepackage{amsrefs}
\usepackage{bigints}

\newtheorem{theorem}{Theorem}[section]

\newtheorem{proposition}[theorem]{Proposition}
\newtheorem{corollary}[theorem]{Corollary}
\theoremstyle{definition}
\newtheorem{definition}[theorem]{Definition}

\theoremstyle{remark}
\newtheorem{remark}[theorem]{Remark}
\numberwithin{equation}{section}

\allowdisplaybreaks

\begin{document}
	\setcounter{page}{1}
	
	\title[]{Sharp Strichartz estimates for some variable coefficient Schr\"{o}dinger operators on $\mathbb{R}\times\mathbb{T}^2$}

	\author[S. Federico]{Serena Federico}
	\address{
		Serena Federico:
		\endgraf
		Department of Mathematics: Analysis Logic and Discrete
		\endgraf
		Mathematics
		\endgraf
		Ghent University
		\endgraf
		Krijgslaan 281, Ghent, B 9000
		\endgraf
		Belgium
		\endgraf
		{\it E-mail address} {\rm serena.federico@ugent.be}}
	\thanks{S.F. has received funding from the European Unions Horizon 2020 research and innovation programme under the Marie Sk\l{}odowska-Curie grant agreement No 838661 and by the FWO Odysseus 1 grant G.0H94.18N: Analysis and Partial Differential Equations. }

	\author[G. Staffilani]{Gigliola Staffilani}
	\address{
		Gigliola Staffilani:
		\endgraf
		Department of Mathematics
		\endgraf
Massachusetts Institute of Technology
		\endgraf
			77 Massachusetts Ave,
		MA 02139-4307
		\endgraf
		USA
		\endgraf
		{\it E-mail address} {\rm  gigliola@math.mit.edu}}
	\endgraf
	\thanks{ G.S. is  funded in part by  DMS-1764403 and by the Simons Foundation}

	\begin{abstract} In the first part of the  paper we continue   the study of solutions to Schr\"odinger equations with a time singularity in the dispersive relation and in the periodic setting. In the second we show that if the Schr\"odinger operator involves a Laplace  operator with variable coefficients with a particular dependence on the space variables, then one can prove Strichartz estimates at the same regularity as that  needed for  constant coefficients.  	 Our work presents a two dimensional analysis, but we expect that with the obvious adjustments similar results are available in higher dimensions.	
	\end{abstract}

\maketitle

\tableofcontents

\section{Introduction} \label{Intro}

This paper is devoted to the analysis of Schr\"{o}dinger operators on the two dimensional torus $\mathbb{T}^2$ with variable coefficients depending either on time or on space.

In the time-dependent case, being, in particular, a time-degenerate case, a time-weighted version of Bourgain's sharp $L^4$-Strichartz estimate will be derived and applied to the corresponding semilinear	initial value problem (IVP).
As for the (special) space-dependent nondegenerate case, the sharp local well-posedness in $H^\varepsilon(\mathbb{T}^2)$, for any $\varepsilon>0$,  of the corresponding nonlinear cubic  initial value problem will be shown.

 Before describing in detail the operators under consideration and the problems addressed in this work, we recall that the celebrated sharp $L^4$-Strichartz estimate on  $\mathbb{T}^2$ proved by Bourgain  \cite{Bourgain-Geom.Funct.Anal-1993} for the rational torus and by Bourgain-Demeter \cite{Bourgain-Demeter} in the general case states that:
\begin{equation}\label{Bou.estimate}
	\| e^{it\Delta_x}u_0\|_{L^4_t([0,1]\times \mathbb{T}^2)}\leq C_\epsilon \| u_0\|_{H^\epsilon(\mathbb{T}^2)},\quad \varepsilon >0.
\end{equation}
Here and in the rest of the paper we shall refer to \eqref{Bou.estimate} as {\it Bourgain's sharp Strichartz} estimate on $\mathbb{T}^2$.

Let us start by introducing Schr\"{o}dinger equations  with time variable coefficients, a  topic that has  attracted  the interest of mathematicians and physicists in the last decades.
For instance, equations of the form
\begin{equation}\label{eq.intro}
i\partial_tu +b(t)\Delta_xu+h(t)u|u|^2+v(t,x)u=0, \
\end{equation}
under some assumptions on the time dependent functions, have been recently studied in the context of Bose-Einstein condensations and nonlinear optics in the Euclidian setting. Exact one and two soliton solutions for \eqref{eq.intro} have been derived in \cite{Lie_2011}, whereas similar related problems have been studied, for example, in \cite{Zheng-Chun-Long}, \cite{Wang_2011} and \cite{ChenYuan-Ming} (see also references therein). Interestingly  it was shown in  \cite{Zheng-Chun-Long} that the nature of the soliton solution (canonical soliton or deformed canonical soliton) to equation \eqref{eq.intro} is influenced by the choice of the time dependent coefficients $b, h$ and $v$.

Time-degenerate equations of the form \eqref{eq.intro} (with $b(0)=0$), still in the Euclidean setting, have been investigated in connection with other mathematical problems. 
Local well-posedness results for the homogeneous IVP associated with time-degenerate equations of the form  \eqref{eq.intro} (with $h=v=0$)  have been derived in \cite{Cicognani-Reissig}. The validity of smoothing estimates and the local well-posedness of the nonlinear IVP for time-degenerate equations of the form  \eqref{eq.intro} was proved by the authors in \cite{Federico-Staffilani}, while Strichartz estimates for the same kind of time-degenerate equations (with $v=0$) have been derived in \cite{Federico-Ruzhansky} by the first author and Ruzhansky. All the results listed hold for space variable in $\mathbb{R}^d$. In the present work instead we are interested on investigating Strichartz estimates and well-posedness in the periodic setting.

Let us now go back to the {\it Bourgain's sharp Strichartz} estimate \eqref{Bou.estimate}. We would like to stress that, in spite of the fact that  a deep analysis of Strichartz estimates on the torus and on more general manifolds has been carried out in the last years, as we will recall below, Bourgain's sharp result remains unproved for variable coefficient Schr\"{o}dinger operators on $\mathbb{T}^2$, even in the simplest case when the Laplace operator is replaced by an elliptic operator with coefficients which are smooth and almost constant.

In the Euclidean setting  Strichartz estimates  for non-degenerate space-variable coefficient  Schr\"{o}dinger operators were proved by the second author and Tataru in \cite{Staffilani-Tatar_2002}, and, in the non-elliptic case, by Salort in \cite{Salort_2007}. Smoothing estimates for non-degenerate space-variable coefficient  Schr\"{o}dinger operators have been considered by Kenig, Ponce, Rolvung and Vega, see, for instance, \cite{Kenig-Ponce-Rolvung-Vega_AdvMath_2005} and references therein.
Smoothing and Strichartz estimates for asymptotically flat Schr\"{o}dinger operators were derived by Marzuola, Metcalfe and Tataru in \cite{Marzuola-Metcalfe-Tataru_2008}.

In the manifolds setting Strichartz estimates on general compact manifolds have been proved by Burq, Gerard and Tzvetkov in \cite{Burq-Gerard-Tvetkov_2004} and by Hani in \cite{Hani}, and by Mizutani and Tzvetkov in \cite{Mizutani-Tzvetkov}  in the non-elliptic case.
Specifically, Theorem 1 in \cite{Burq-Gerard-Tvetkov_2004} states the following:
\medskip

{\it
 Let $(M,g)$ be a Riemannian compact manifold of dimension $d\geq1$ and $\Delta$  be the Laplace Beltrami operator on $M$. Given $p, q$ satisfying the scaling condition $2/p+d/q=d/2$, and $p\geq 2$, $q<\infty$, the solution $v$ of 
	$$i\partial_t v(t,x)+\Delta v(t,x)=0,\quad v(0,x)=v_0(x),$$
satisfies for any finite time interval $I$,
\begin{equation}\label{strichartz.BGT}
\|v(t)\|_{L^p(I)L^q(M)}\leq C(I)\| v_0\|_{H^{\frac 1 p}(M)}.
\end{equation}
}

It is clear that Boungain's sharp result on $\mathbb{T}^2$ is not covered by \eqref{strichartz.BGT}, since, for $p=q=2$ and $M=\mathbb{T}^2$, \eqref{strichartz.BGT} gives
$$\|v(t)\|_{L^4(I)L^4(\mathbb{T}^2)}\leq C(I)\| v_0\|_{H^{\frac 1 4}(\mathbb{T}^2)}.$$
However, still in \cite{Burq-Gerard-Tvetkov_2004}, it is shown that estimate \eqref{strichartz.BGT} is actually sharp when $M$ is the sphere $\mathbb{S}^d$ with $d\geq 3$.

The  observations and results mentioned above motivated our analysis of the following linear  Schr\"{o}dinger equations with variable coefficients  in  $\mathbb{T}^2$:
\begin{equation}\label{time.S}
i\partial_t u+g'(t)\Delta_x u=0,
\end{equation}
and
\begin{equation}\label{space.S}
i\partial_t u+a_1(x_1)\partial_{x_1}^2 u+a_2(x_2)\partial_{x_2}^2 u=0,
\end{equation}
where $g'$, representing the derivative of a strictly monotone function $g\in C^1(\mathbb{R})$, with $g(0)=0$, is such that $g'(0)=0$, while $a_i\in C^\infty(\mathbb{T})$, for $i=1,2$, are strictly positive functions\footnote{Here we assume smoothness for $a_1$ and $a_2$, but much less regularity is required for our analysis.}. 

For the time-degenerate equation \eqref{time.S} (which is now considered in a non-Euclidean setting),  a suitable {\it weighted} version of Bourgain's sharp-Strichartz estimate on $\mathbb{T}^2$ will be obtained.
As a consequence of the aforementioned inequality, a local well-posedness result for the semilinear IVP associated with \eqref{time.S} will be obtained.

As for  equation \eqref{space.S}, by exploiting a change of variables in combination with the use of a gauge transform, we will be able to apply Bourgain's sharp-Strichartz estimate to obtain the local well-posedness of the corresponding cubic IVP in $H^\varepsilon(\mathbb{T}^2), \, \varepsilon>0$ (the same result as the one we observe in the non-degenerate constant coefficients case). 

Although  the study of \eqref{time.S} and \eqref{space.S} is interesting on its own term, the analysis of these operators show that sharp Strichartz estimates on $\mathbb{T}^2$ can be obtained even when some variable (possibly degenerate in time) coefficients are present, and as a consequence also  sharp local well-posedness results in the same settings.

Indeed, while the result in \cite{Burq-Gerard-Tvetkov_2004}
allows to prove local well-posedness in $H^s(\mathbb{T}^2)$, with $s>1/4$, of the cubic IVP associated with the {\it constant} coefficient Schr\"{o}dinger operator on $\mathbb{T}^2$, here we manage to combine different techniques to obtain sharp local well-posedness in $H^s$, $s >0$, when some variable coefficients are present.  

For completeness we list below the main results of this paper. 

\addtocontents{toc}{\protect\setcounter{tocdepth}{1}}
\subsection*{Well-posedness of cubic Nonlinear Schr\"{o}dinger equations with variable coefficients on $\mathbb{R}\times\mathbb{T}^2$.}

Let $g\in C^\infty(\mathbb{R}_t)$ be a strictly monotone function such that $g(0)=g'(0)=0$, let $f$ be a smooth function such that $f\in H^{1,b}_g(\mathbb{R})$ for $b\in (1/2,1)$, where $H^{1,b}_g(\mathbb{R})$ is as in Definition \ref{def.Hspaces}, and let $\tilde{X}^{s,b}_g(\mathbb{R}\times \mathbb{T}^2)$ be the Banach space as in Definition \ref{def.Xsbg}. Then the following theorems hold. 

\begin{theorem}\label{thm.time-lwp}
	Let $s>0$ and $b\in (1/2,1)$. Then, for every $u_0\in H^s(\mathbb{T}^2)$, there exists  a unique solution of the IVP 
	\begin{equation}\label{IVPtimeIntro}
	\left\{
	\begin{array}{l}
	i\partial_t u + g'(t)\Delta_x u= g'(t)|u|^2u,\\
	u(0,x)=u_0(x),
	\end{array}
	\right.
	\end{equation}
	in the time interval $[-T,T]$ for a suitable time $T=T(\|u_0\|_{H^s})$.  Moreover the solution $u$ satisfies
	$$u\in C([-T,T];H^s)$$
	and, for $I$ closed neighborhood of $[-T,T]$, and $\chi_I$ a smooth cutoff function such that $\chi_I\equiv 1$ on $[-T,T]$, we have that there exists $b\in (1/2,1)$ such that
	
	$$\chi_I u\in \tilde{X}^{s,b}_g(\mathbb{R}\times \mathbb{T}^2).$$
\end{theorem}

\begin{theorem}\label{thm.time-lwp2}
	Let $s>0$ and $b\in (1/2,1)$. Then, for every $u_0\in H^s(\mathbb{T}^2)$, there exists  a unique solution of the IVP  
	\begin{equation}\label{IVPtime2Intro}
	\left\{
	\begin{array}{l}
	i\partial_t u + g'(t)\Delta_x u= f(t)|u|^2u,\\
	u(0,x)=u_0(x),
	\end{array}
	\right.
	\end{equation}
	in the time interval $[-T,T]$ for a suitable time $T=T(\|u_0\|_{H^s})$.  Moreover the solution $u$ satisfies
	$$u\in C([-T,T];H^s)$$
	and, for a closed neighborhood $I$ of $[-T,T]$, we have that there exists $b\in (1/2,1)$ such that
	$$\chi_I u\in \tilde{X}^{s,b}_g(\mathbb{R}\times \mathbb{T}^2)$$
	with $\chi_I$ being  a smooth cutoff function such that $\chi_I\equiv 1$ on $[-T,T]$.
\end{theorem}

Let $a_1,a_2\in C^\infty (\mathbb{T})$ be two real valued strictly positive functions and $X^{s,b}_{\Phi,\tilde{ \alpha}}(\mathbb{R}\times \mathbb{T}^2)$ the Banach space in Definition \ref{def.PhiXsbg}. 
Then the following local well-posedness result holds.
\begin{theorem}\label{thm.lwp.space-Intro}
Let $s>0$ and $b\in (1/2,1)$. Then, for every $u_0\in H^s(\mathbb{T}^2)$, there exists  a unique solution of the IVP  
\begin{equation}\label{IVP3Intro}
\left\{	\begin{array}{l}
i\partial_t u+a_1(x_1)\partial_{x_1}^2 u+a_2(x_2)\partial_{x_2}^2 u=  u|u|^2,\\
u(0,x)=u_0(x),
\end{array}\right.
\end{equation}
	in the time interval $[-T,T]$ for a suitable time $T=T(\|u_0\|_{H^s})$.  Moreover the solution $u$ satisfies
$$u\in C([-T,T];H^s)$$
and, for a closed neighborhood $I$ of $[-T,T]$, we have  that there exists $b\in (1/2,1)$ such that
$$\chi_I u\in X_\Phi^{s,b}(\mathbb{R}\times \mathbb{T}^2)$$
with $\chi_I$ being  a smooth cutoff function such that $\chi_I\equiv 1$ on $[-T,T]$.
\end{theorem}

\begin{remark}
	Applying the strategies employed in the results above one can also derive the local well-posedness of the IVP
\begin{equation*}
\left\{	\begin{array}{l}
i\partial_t u+g'(t)\Big(a_1(x_1)\partial_{x_1}^2 u+a_2(x_2)\partial_{x_2}^2 u\Big)=  f(t)u|u|^2,\\
u(0,x)=u_0(x),
\end{array}\right.
\end{equation*}
with solution in a suitable   $X^{s,b}_{g,\Phi}$ space (see Definition \ref{def.PhiXsbg}).
\end{remark}

\subsection*{Well-posedness of the quintic Nonlinear Schr\"{o}dinger equations on $\mathbb{R}\times\mathbb{T}$.}
Under the same assumptions on the functions $g$ and $f$ the following results hold.

\begin{theorem}\label{thm.time-lwp-Q}
	Let $s>0$ and $b\in (1/2,1)$. Then, for every $u_0\in H^s(\mathbb{T})$, there exists  a unique solution of the IVP 
	\begin{equation}\label{IVPtimeIntro-Q}
	\left\{
	\begin{array}{l}
	i\partial_t u + g'(t)\Delta_x u= g'(t)|u|^4u,\\
	u(0,x)=u_0(x),
	\end{array}
	\right.
	\end{equation}
	in the time interval $[-T,T]$ for a suitable time $T=T(\|u_0\|_{H^s})$.  Moreover the solution $u$ satisfies
	$$u\in C([-T,T];H^s)$$
	and, for $I$ closed neighborhood of $[-T,T]$, and $\chi_I$ a smooth cutoff function such that $\chi_I\equiv 1$ on $[-T,T]$, we have that there exists $b\in (1/2,1)$ such that
	
	$$\chi_I u\in \tilde{X}^{s,b}_g(\mathbb{R}\times \mathbb{T}).$$
\end{theorem}

\begin{theorem}\label{thm.time-lwp2-Q}
	Let $s>0$ and $b\in (1/2,1)$. Then, for every $u_0\in H^s(\mathbb{T})$, there exists  a unique solution of the IVP  
	\begin{equation}\label{IVPtime2Intro-Q}
	\left\{
	\begin{array}{l}
	i\partial_t u + g'(t)\Delta_x u= f(t)|u|^4u,\\
	u(0,x)=u_0(x),
	\end{array}
	\right.
	\end{equation}
	in the time interval $[-T,T]$ for a suitable time $T=T(\|u_0\|_{H^s})$.  Moreover the solution $u$ satisfies
	$$u\in C([-T,T];H^s)$$
	and, for a closed neighborhood $I$ of $[-T,T]$, we have that there exists $b\in (1/2,1)$ such that
	$$\chi_I u\in \tilde{X}^{s,b}_g(\mathbb{R}\times \mathbb{T})$$
	with $\chi_I$ being  a smooth cutoff function such that $\chi_I\equiv 1$ on $[-T,T]$.
\end{theorem}

Let $a\in C^\infty (\mathbb{T})$ be a real valued strictly positive function and $X^{s,b}_{\Phi,\tilde{ \alpha}}(\mathbb{R}\times \mathbb{T})$ as in Definition \ref{def.PhiXsbg}. 
Then the following local well-posedness result holds.
\begin{theorem}\label{thm.lwp.space-IntroQ}
	Let $s>0$ and $b\in (1/2,1)$. Then, for every $u_0\in H^s(\mathbb{T})$, there exists  a unique solution of the IVP  
	\begin{equation}\label{IVP3Intro-Q}
	\left\{	\begin{array}{l}
	i\partial_t u+a(x)\Delta_x u=  u|u|^4,\\
	u(0,x)=u_0(x),
	\end{array}\right.
	\end{equation}
	in the time interval $[-T,T]$ for a suitable time $T=T(\|u_0\|_{H^s})$.  Moreover the solution $u$ satisfies
	$$u\in C([-T,T];H^s)$$
	and, for a closed neighborhood $I$ of $[-T,T]$, we have that there exists $b\in (1/2,1)$ such that
	$$\chi_I u\in X_{\Phi,\tilde{ \alpha}}^{s,b}(\mathbb{R}\times \mathbb{T})$$
	with $\chi_I$ being  a smooth cutoff function such that $\chi_I\equiv 1$ on $[-T,T]$.
\end{theorem}
\begin{remark}
	Applying the strategies employed in the results above one can also derive the local well-posedness of the IVP
	\begin{equation*}
	\left\{	\begin{array}{l}
	i\partial_t u+g'(t)a(x)\Delta_x u=  f(t)u|u|^4,\\
	u(0,x)=u_0(x),
	\end{array}\right.
	\end{equation*}
	with solution in a suitable   $\tilde{X}^{s,b}_{g,\Phi,\tilde{ \alpha}}$ space (see Definition \ref{def.PhiXsbg}).
\end{remark}
\medskip

Below we will give detailed proofs of the two dimensional case ($x\in\mathbb{T}^2$), that  have straightforward applications in the one dimensional case ($x\in\mathbb{T}$). However we expect equivalent  results to be true in higher dimensions after applying suitable adjustments in our assumptions and arguments.

Let us remark that in the sequel we will also give very precise proofs of some standard multilinear estimates contained, for instance, in \cite{Bourgain-Geom.Funct.Anal-1993} and \cite{Kenig-Ponce-Rolvung-Vega_AdvMath_2005}.  Even if some of these results are well-known and standard, detailed proofs allow to measure in a more precise way some key parameters appearing in the contraction argument, as, for instance, the parameter $\delta$ associated with the length of the  time interval of existence, and, more importantly, some (dangerous) exponents appearing in the contraction argument. It is only with a sharp evaluation of these quantities that the contraction argument can be performed successfully and get {\it sharp} local well-posedness results.

\medskip

We conclude this introduction by giving the plan of the paper. Section \ref{sec.prelim} is devoted to the definition of the $X^{s,b}_g$-spaces which are central in our analysis.
In Section \ref{sec.time} we derive the suitable weighted formulation of Strichartz and multilinear estimates via de use of $X^{s,b}_g$-spaces. This section also contains the proof of Theorem \ref{thm.time-lwp}, \ref{thm.time-lwp2}, \ref{thm.time-lwp-Q} and \ref{thm.time-lwp2-Q}.
Finally Section \ref{Section.space.variable} will be devoted to the study of \eqref{IVP3Intro} and \eqref{IVP3Intro-Q}, that is, in particular, to the proof of Theorem \ref{thm.lwp.space-Intro} and \ref{thm.lwp.space-IntroQ}.

\section{Preliminaries}\label{sec.prelim}
In this section we shall briefly recall some tools we will be using throughout the paper. 

Given a function $u\in  C^\infty(\mathbb{T}^d)$, the Fourier and anti-Fourier transforms of $u$ are defined as
$$\widehat{u}(k)=(\mathcal{F}u)(k)=\int_{\mathbb{T}^d} e^{-i k\cdot x}u(x)dx,$$
and
$$u(x)=(\mathcal{F}^{-1}\widehat{u})(x)=\sum_{k\in\mathbb{Z}^{d}}e^{i k\cdot x}\widehat{u}(k),$$
and, by Plancherel's Theorem, we have the identity
$$\| u\|_{L^2(\mathbb{T}^d)}=\| \widehat{u} \|_{\ell^2(\mathbb{Z}^d)}.$$

Similarly, given a space-time dependent function $v(t,x)$, with $(t,x)\in \mathbb{R}\times \mathbb{T}^d$,  we can write $v(t,x)$ by means of the Fourier representation formula as
$$v(t,x)= \sum_{k\in\mathbb{Z}^{d}}\int_{\mathbb{R}}e^{-i(k\cdot x +t\tau)}\widehat{v}(k,\tau)d\tau,$$
where here $\widehat{v}(k,\tau)$ (and in the rest of the paper) is the space-time Fourier transform.

As Bourgain did in his pioneering work in the periodic setting \cite{Bourgain-Geom.Funct.Anal-1993}, also in our case we will  make use of the $X^{s,b}$-spaces  recalled below.
\medskip

\begin{definition}[$X^{s,b}$ Spaces] \label{def.Xsb}
Let $X$ be the space of functions $u$ on $\mathbb{R}\times \mathbb{T}^d$ such that

\begin{itemize}
	\item [(i)] $u:\mathbb{R}\times \mathbb{T}^d \rightarrow \mathbb{C}$;
	\item[(ii)] $ u(\cdot,x)\in \mathcal{S}(\mathbb{R})$ for all $x\in \mathbb{T}^d$;
	\item[(iii)] $x\rightarrow u(\cdot,x)$ is $C^\infty(\mathbb{T}^d)$.
\end{itemize}
Then, for $s,b\in \mathbb{R}$, we denote by $X^{s,b}$ the completion of the space $X$ with respect to the norm
$$\|u\|_{X^{s,b}}:=\left(\sum_{k\in\mathbb{Z}^d}(1+|k|)^{2s}\int_{\mathbb{R}}(1-|\tau-|k|^2|)^{2b}|\widehat{u}(\tau,k) |^2d\tau \right)^{1/2}.$$
\end{definition}

We will now introduce a new notion of $X^{s,b}$ space subordinate to a function $g\in C^\infty(\mathbb{R}_t)$. 
Let $g\in C^\infty(\mathbb{R}_t)$ be a strictly monotone function such that $g(0)=0$, then we define the  modified Fourier and inverse Fourier transform subordinate to $g$ as
$$(\tilde{\mathcal{F}}u)(\tau):=\int_\mathbb{R}e^{-ig(t)\tau}u(t) dt$$
and 
$$(\tilde{\mathcal{F}}^{-1}v)(t):=g'(t)\int_\mathbb{R}e^{ig(t)\tau}v(\tau) d\tau,$$
where $u,v\in L^1(\mathbb{R})\cap L^2(\mathbb{R})$.
Note that, on denoting by $\tilde{u}(\tau):=(\tilde{\mathcal{F}}u)(\tau)$ we have that the following properties hold
\begin{itemize}
\item $u(t)=(\tilde{\mathcal{F}}^{-1}\tilde{u})(t); $
\item $\tilde{\mathcal{F}}(\partial_t u)(\tau)=(-i\tau)\tilde{\mathcal{F}}(g'u)(\tau)$;
\item $\|\tilde{u}\|_{L^2(\mathbb{R}_\tau)}=\|\frac{1}{\sqrt{|g'|}}u\|_{L^2(\mathbb{R}_t)}$;
\item $\|\widetilde{g'u}\|_{L^2(\mathbb{R}_\tau)}=\|\sqrt{|g'|}\,u\|_{L^2(\mathbb{R}_t)}$;
\item $\|\tilde{\mathcal{F}}(\sqrt{|g'|}u)\|_{L^2(\mathbb{R}_\tau)}=\|u\|_{L^2(\mathbb{R}_t)}$;
\end{itemize}
When no confusion arises we shall also use the notation $\tilde{u}(\tau,k)$ for the space-time transform
$$\tilde{u}(\tau,k):=\int_{\mathbb{R}\times\mathbb{T}^d} e^{-i(g(t)\tau+k\cdot x)}u(t,x)dt dx,$$
namely, the modified Fourier transform in time and the standard Fourier transform in space of a function $u$ on $\mathbb{R}\times \mathbb{T}^d$.

With this definition at our disposal we can now define the $X^{s,b}_g$ and the $\tilde{X}^{s,b}_g$ spaces subordinate to $g$.
\medskip

\begin{definition}[$X^{s,b}_g$ and $\tilde{X}^{s,b}_g$ spaces]\label{def.Xsbg}
Given a strictly monotone function $g\in C^\infty(\mathbb{R}_t)$, with $g(0)=0$, we define the space  $X^{s,b}_g$ as the completion of the space $X$ in Definition \ref{def.Xsb} with respect to the norm
$$\|u\|_{X^{s,b}_g}:=\left(\sum_{k\in\mathbb{Z}^d}(1+|k|)^{2s}\int_{\mathbb{R}}(1+|\tau-|k|^2|)^{2b}|\tilde{u}(\tau,k) |^2d\tau \right)^{1/2}.$$
With the same function $g$ we define the spaces $\tilde{X}^{s,b}_g$ as 

$$\tilde{X}^{s,b}_g:=\{u\in X^{s,b}_g; g'u\in X^{s,b}_g\}.$$
\end{definition}

\begin{remark}
Note that, by a simple change of variables, we have the following relations between $X^{s,b}_g$, $\tilde X^{s,b}_g$ and $X^{s,b}$ spaces:
$$\| u\|_{X^{s,b}_g}=\left\| \frac{1}{g'\circ g^{-1}}\,\, u(g^{-1}(\cdot),\cdot)\right \|_{X^{s,b}}$$
and
$$\|u\|_{\tilde X^{s,b}_g}:=\|g'u\|_{X^{s,b}_g}=\|u\left(g^{-1}(\cdot),\cdot\right)\|_{X^{s,b}},$$
provided that $g'u$ and $u$ are smooth enough.
Note also that
$$\| u\|_{X^{0,0}_g}=\|\tilde{u}\|_{L^2_\tau\ell^2_k}=\left\|\frac{1}{\sqrt{|g'|}}u\right \|_{L^2_{t,x}}.$$

These simple observations will be crucial in the next section where the time-degenerate case described by equation \eqref{time.S} will be treated. We shall also explain below how these spaces are related with the solution of the IVP associated with \eqref{time.S}.
\end{remark}

We conclude this section by defining some other spaces we will be using throughout the paper.

\begin{definition}[$X^{s,b}_\Phi$, $X^{s,b}_{g,\Phi}$, $X^{s,b}_{\Phi,\tilde{ \alpha}}$, $X^{s,b}_{g,\Phi,\tilde{\alpha}}$ , and $\tilde{X}^{s,b}_{g,\Phi,\tilde{ \alpha}}$ spaces]
\label{def.PhiXsbg}
Let $\Phi \in C^\infty(\mathbb{T}^d)$ and let $g$ be as in Definition \ref{def.Xsbg}. Let also $\tilde{\alpha}: \mathbb{R}\times \mathbb{T}^d\rightarrow \mathbb{R}\times \mathbb{T}^d$ be such that $\tilde{\alpha}(t,x):=(t, \alpha(x))$, where $\alpha:\mathbb{T}^d\rightarrow \mathbb{T}^d$ is a diffeomorphism.
Then we define the spaces $X^{s,b}_\Phi$, $X^{s,b}_{g,\Phi}$, $X^{s,b}_{\Phi,\tilde{ \alpha}}$, $X^{s,b}_{g,\Phi,\tilde{\alpha}}$ , and $\tilde{X}^{s,b}_{g,\Phi,\tilde{ \alpha}}$ as
\begin{eqnarray*}
X^{s,b}_\Phi(\mathbb{R}\times \mathbb{T}^d)&:=&\{f:\mathbb{R}\times \mathbb{T}^d\rightarrow \mathbb{C}; e^{\Phi}f\in X^{s,b}(\mathbb{R}\times \mathbb{T}^d)\},\\
X^{s,b}_{g,\Phi}(\mathbb{R}\times \mathbb{T}^d)&:=&\{f:\mathbb{R}\times \mathbb{T}^d\rightarrow \mathbb{C}; e^{\Phi}f\in X^{s,b}_g(\mathbb{R}\times \mathbb{T}^d)\},\\
\tilde{X}^{s,b}_{g,\Phi}(\mathbb{R}\times \mathbb{T}^d)&:=&\{f:\mathbb{R}\times \mathbb{T}^d\rightarrow \mathbb{C}; g'e^{\Phi}f\in X^{s,b}_g(\mathbb{R}\times \mathbb{T}^d)\},\\
X^{s,b}_{\Phi,\tilde{\alpha}}(\mathbb{R}\times \mathbb{T}^d)&:=&\{f:\mathbb{R}\times \mathbb{T}^d\rightarrow \mathbb{C}; (e^{\Phi}\,f)\circ \tilde{\alpha}\in X^{s,b}(\mathbb{R}\times \mathbb{T}^d)\},\\
X^{s,b}_{g,\Phi,\tilde{\alpha}}(\mathbb{R}\times \mathbb{T}^d)&:=&\{f:\mathbb{R}\times \mathbb{T}^d\rightarrow \mathbb{C}; (e^{\Phi}\,f)\circ \tilde{\alpha}\in X^{s,b}_g(\mathbb{R}\times \mathbb{T}^d)\},\\
\tilde{X}^{s,b}_{g,\Phi,\tilde{\alpha}}(\mathbb{R}\times \mathbb{T}^d)&:=&\{f:\mathbb{R}\times \mathbb{T}^d\rightarrow \mathbb{C}; (e^{\Phi}\,f)\circ \tilde{\alpha}\in \tilde{X}^{s,b}_g(\mathbb{R}\times \mathbb{T}^d)\}.
\end{eqnarray*}
\end{definition}

\begin{definition}[$H^{p,b}$ and $H^{p,b}_g$ spaces]\label{def.Hspaces}
	Let $p\in [1,\infty)$ and $b\in \mathbb{R}$, then we define the spaces $H^{p,b}(\mathbb{R})$ and $H^{p,b}_g(\mathbb{R})$ as
	$$H^{p,b}(\mathbb{R}):=\{f\in L^p(\mathbb{R}); \widehat{f}, \widehat{D^bf}\in L^p(\mathbb{R})\}$$
	equipped with the norm
	$$\| f\|^p_{H^{p,b}}:=\int_{\mathbb{R}} \langle \tau\rangle^{pb} |\widehat{f}(\tau)|^p d\tau,$$ 
	with $\langle \tau \rangle :=(1+|\tau|^2)^{1/2}$, and 
	$$H^{p,b}_g(\mathbb{R}):=\{f\in L^p(\mathbb{R}); \| f\|_{H^{p,b}_g}<\infty \},$$ where $\| f\|^p_{H^{p,b}_g}:=\int_{\mathbb{R}} \langle \tau\rangle^{pb} |\tilde{f}(\tau)|^p d\tau$. 
\end{definition}

\section{The degenerate time-dependent  case}\label{sec.time}
Here we focus on the analysis of the two-dimensional case, that is when $x\in \mathbb{T}^2$.  We will use the square torus  to conduct our calculations, so that the symbol of the Laplacian is  simply $-|k|^2$. The argument will be the same for any other torus.  The one-dimensional case follows similarly and is described at the end of this section.

We consider the time-degenerate Schr\"{o}dinger equation on $\mathbb{R}\times \mathbb{T}^2$ 
$$i\partial_t u-g'(t)\Delta_x u=0,$$
where $g'$ is the derivative of a strictly monotone function $g\in C^{\infty}(\mathbb{R})$, with $g(0)=0$, and such that $g'(0)=0$. Additionally, $g$ is supposed to have power growth, that is there exists $\alpha> 0
$ such that $|g^{(j)}(t)|\lesssim (1+|t|)^{\alpha-j}$.

It is easy to see that the solution of the IVP associated with the equation under consideration and with initial datum $u(0,x)=u_0$ is given by

$$u(t,x)=S(t)u_0:=\sum_{k\in\mathbb{Z}^{2}}e^{-ig(t)|k|^2}\widehat{u}_0(k),$$
where $S(t):=S(t,0):=e^{ib(t)\Delta_x}$ represents the so called {\it solution operator } giving the solution at time $t$ starting at time $0$. More generally, given a space-dependent function $\varphi$, we have

$$S(t,s)\varphi :=\sum_{k\in\mathbb{Z}^{2}}e^{-i(g(t)-g(s))|k|^2}\widehat{\varphi}(k),$$
where
\begin{itemize}
	\item $S(t,t)=I\quad \forall t\in \mathbb{R}$;
	\item $S(t,s)=S(t,r)S(r,s),\quad \forall r,s,t\in \mathbb{R}$;
	\item $S(t,s)\Delta_x=\Delta_x S(t,s).$
\end{itemize}
Moreover by Duhamel's principle  the solution of the inhomogeneous IVP
$$
\left\{ \begin{array}{l}
i\partial_t u-g'(t)\Delta_x u=f(t,x) \\
u(0,x)=u_0(x)
\end{array}
\right.
$$
is given by
$$\tilde{S}(t)u_0:=S(t)u_0+\int_{0}^{t} S(t,s)f(s)ds.$$

Now, considering $\tilde{\mathcal{F}}$ as the time-Fuorier transform  subordinate to $g$ and applying $\tilde{\mathcal{F}}_{t\rightarrow\tau}\mathcal{F}_{x\rightarrow k}$ to equation \eqref{time.S}, we obtain that 
$$\tilde{\mathcal{F}}_{t\rightarrow\tau}\mathcal{F}_{x\rightarrow k}(i\partial_t u-g'(t)\Delta_x u)=0$$
which is equivalent to
$$(-\tau +|k|^2)\widetilde{g' u}(\tau,k)=0,$$
meaning that the modified Fourier transform of the product of $g'$ and the solution $u$ is supported on the paraboloid $\tau=|k|^2$.

\addtocontents{toc}{\protect\setcounter{tocdepth}{1}}
\subsection{Strichartz and multilinear estimates}
\addtocontents{toc}{\protect\setcounter{tocdepth}{2}}
Under the previous assumptions on the function $g$ we now prove the frequency localized (time-)weighted Strichartz estimate from which the full (time-)weighted Strichartz estimate will follow. Afterwards, we shall also give detailed proofs of some classical multilinear estimates holding in the constant coefficients case. We will then translate these estimates in their suitable (time-)weighted version to be used in the time-degenerate case.
Some of the standard multilinear estimates prove here will be used in their original form in the analyis of the space-variable coefficients case studied in Section \ref{Section.space.variable}.

For simplicity we assume that $g$ is strictly increasing, where, recall, $g(0)=g'(0)=0$. This assumption is useful when we make certain change of variables.
Note that the vanishing of $g'$ is essential in order to have a degeneracy in the time-dependent Schr\"{o}dinger operator.
\medskip

\begin{proposition}
Let $\phi_N\in L^2(\mathbb{T}^2)$ be such that $\mathrm{supp}\, \widehat{\phi}_N\subseteq B(0,N):=\{k\in\mathbb{Z}^2; |k|\leq N\}$, and let $I\subset\mathbb{R}$ be a finite interval centered at zero.
Then, for any $\epsilon>0$, we have
\begin{equation}\label{loc.strich}
\| g'(t)^{1/4}S(t)\phi_N\|_{L^4(I\times \mathbb{T}^2)}
\lesssim  N^\epsilon\|\phi_N\|_{L^2(\mathbb{T}^2).}
\end{equation}
\end{proposition}
\proof
To prove the result we apply a change of variable in time and get

$$\| g'(t)^{1/4}S(t)\phi_N\|_{L^4(I\times \mathbb{T}^2)}
=\left(\int_{I\times\mathbb{T}^d}  g'(t)|e^{ig(t)\Delta}\phi_N(x)|^4 dx dt\right)^{1/4}$$
$$\underset{t'=g(t)}{=}\left(\int_{I'\times\mathbb{T}^d}    |e^{it'\Delta}\phi_N(x)|^4 dx dt'\right)^{1/4}=\|  e^{it'\Delta}\phi_N\|_{L^4(I'\times \mathbb{T}^2)}$$
$$\lesssim  N^\epsilon \|\phi_N\|_{L^2(\mathbb{T}^2)},$$
where the last inequality follows from the standard case when $g'(t)=1$ (see \cite{Bourgain-Geom.Funct.Anal-1993} for the rational case and \cite{Bourgain-Demeter} for the general case). This shows \eqref{loc.strich} and concludes the proof.
\endproof

\begin{proposition}\label{thm.loc.Str}
Let $\phi_{N_1},\phi_{N_2}\in L^2(\mathbb{T}^2)$, $N_1>>N_2$, such that $\mathrm{supp}\, \widehat{\phi}_{N_i}\subseteq B(0,N_i):=\{k\in\mathbb{Z}^2; |k|\leq N_i\}$, and let $I\subset\mathbb{R}$ be a finite interval centered at zero.
Then, for any $\varepsilon>0$, 
\begin{equation}\label{bilinear.est}
\|  g'(t)^{1/2}(S(t)\phi_{N_1} )(S(t)\phi_{N_2})\|_{L^2(I\times \mathbb{T}^2)}
\lesssim  \min(N_1,N_2)^\epsilon\|\phi_{N_1}\|_{L^2(\mathbb{T}^2)}\|\phi_{N_2}\|_{L^2(\mathbb{T}^2)}.
\end{equation}.
\end{proposition}
\proof
The proof is straightforward and follows by simply observing that
$$\| g'(t)^{1/2}(S(t)\phi_{N_1} )(S(t)\phi_{N_2})\|_{L^2(I\times \mathbb{T}^2)}=
\| (e^{it\Delta}\phi_{N_1})(e^{it\Delta}\phi_{N_2})\|_{L^2(I'\times \mathbb{T}^2)}$$
$$\lesssim  N_2^\epsilon\|\phi_{N_1}\|_{L^2(\mathbb{T}^2)}\|\phi_{N_2}\|_{L^2(\mathbb{T}^2)},$$
where, once more, the last inequality is due to the application of classical results (see \cite{Bourgain-Geom.Funct.Anal-1993} for rational tori and \cite{Fan-Ou-Staffilani-Wang} for the general case).
\endproof

%

A Strichartz estimate similar to the one in Proposition  \ref{thm.loc.Str} in any dimensions and any $p\geq 2$, with $g(t)=1$, was proved by Bourgain and Demeter \cite{Bourgain-Demeter}, and hence, using the  Littlewood-Paley decomposition one gets the full $L^p$-Strichartz estimate on $I\times \mathbb{T}^2$ for $p\geq 2$.
%
%
%

\begin{theorem}
Let  $I$ be a finite time interval. Then, for $p\geq 2$, 

$$\|  g'(t)^{1/p}S(t)\phi\|_{L^p(I\times\mathbb{T}^d)}
\lesssim \|\phi\|_{L^2(\mathbb{T}^d)}, \quad  p<\frac{2(d+2)}{d};$$

$$\|  g'(t)^{1/p}S(t)\phi\|_{L^p(I\times\mathbb{T}^d)}
\lesssim \|\phi\|_{H^s(\mathbb{T}^d)}, \quad s>0, \,\, p=\frac{2(d+2)}{d};$$

$$\|  g'(t)^{1/p}S(t)\phi\|_{L^p(I\times\mathbb{T}^2)}
\lesssim \|\phi\|_{H^s(\mathbb{T}^2)}, \quad s>\frac{d}{2}-\frac{d+2}{p}, \,\, p>\frac{2(d+2)}{d}.$$
\end{theorem}

We now restrict ourselves again to the case $d=2$. 
\begin{proposition}
Let $u_N$ be a function on $I\times \mathbb{T}^2$ such that the space-Fourier transform $\widehat{u}_N(t,k)$ is supported in $B(0,N):=\{k\in \mathbb{Z}^2; |k|\leq N\}$. Then, for any $s_1>0$ and $b_1>\frac{1- \min\{s_1,1/2\}}{2}$ ($\frac 1 4<b_1< \frac 1 2+$, $s_1>1-2b_1$ ), we have
$$\| \chi_{[0,1]} g'(t)^{1/4} u\|_{L^4(\mathbb{R}\times\mathbb{T}^2)}\lesssim N^{s_1}\| \chi_{[0,1]}g'(t)u \|_{X^{0,b_1}_g}.$$

\end{proposition}

\proof
By applying the previous change of variables we have
$$\| \chi_{[0,1]} g'(t)^{1/4} u\|_{L^4(\mathbb{R}\times\mathbb{T}^2)}=
\| \chi_{[0,g(1)]}  u(g^{-1}(t),\cdot)\|_{L^4(\mathbb{R}\times\mathbb{T}^2)}$$
$$\lesssim N^{s_1} \| \chi_{[0,g(1)]}  u(g^{-1}(t),\cdot)\|_{X^{0,b_1}}=N^{s_1} \| \chi_{[0,1]} g'(t)  u\|_{X^{0,b_1}_g},$$
where in the last line we applied the result in the standard case $g'(t)=1$ in \cite{Bourgain-Geom.Funct.Anal-1993} and the relation between $X^{s,b}$ and $X^{s,b}_g$ norms. 
\endproof

\begin{proposition}\label{prop.xspEst}
Assume that  $|I'|:=|b(I)|=\delta$, then
	\begin{equation} \label{Xsb_est_1}
	\| \chi_I (t) g'(t)S(t) u_0\|_{X^{s,b}_g}\lesssim \delta^{1/2-b}\|u_0\|_{H^s},\quad \forall u_0\in H^s(\mathbb{T}^2), 
		\end{equation}
		
\begin{equation}\label{Xsb_est_2}
	\left	\|  g'(t)\int_0^t g'(s)S(t,s)w(s)ds \right  \|_{X^{s,b}_g}\lesssim \|g'(t) w\|_{X^{s,b-1}_g},
		\end{equation}
		
	\begin{equation}\label{Xsb_est_3}
	\| \chi_I (t) g'(t)  |u|^2u\|_{X^{s,b-1}_g}\lesssim \|\chi_I g'(t) u\|_{X^{s,b'}_g}^2
	\|\chi_I g'(t) u\|_{X^{s,b}_g} ,\quad \text{for}\,\, b>1/2, \,\,1/4<b'<b,\, s>0;
	\end{equation}
	
	\begin{equation}\label{Xsb_est_4}
	\| \chi_I  g'(t) u \|_{X^{s,b'}_g}\lesssim \delta^{\frac{b-b'}{8}}\|g'(t)u\|_{X^{s,b}_g}.
	\end{equation}
\end{proposition}

\proof
First observe that \eqref{Xsb_est_1}, \eqref{Xsb_est_2}, \eqref{Xsb_est_3} and \eqref{Xsb_est_4}    hold true in the standard case when $g(t)=t$ (see \cite{Bourgain-Geom.Funct.Anal-1993}).
Now, by using the definition of modified Fourier transform (in time) with respect to $g$, we have
$$\tilde{\mathcal{F}}_{t\rightarrow \tau}\mathcal{F}_{x\rightarrow k}(\chi_I  g'(t)S(t) u_0)(\tau,k)=\int_\mathbb{R} e^{-ig(t)\tau -ig(t)|k|^2} \chi_I g'(t) \widehat{u}_0(k)dt$$
$$\underset{t'=g(t)}{=} \int_\mathbb{R} e^{-it'(\tau +|k|^2)} \chi_{g(I)} (t')\widehat{u}_0(k)dt'
=\mathcal{F}_{t\rightarrow \tau}\mathcal{F}_{x\rightarrow k}(\chi_{g(I)}e^{it\Delta }u_0)(\tau,k),$$
which, in particular, gives
$$ 	\| \chi_I (t) g'(t)S(t) u_0\|_{X^{s,b}_g}= 	\| \chi_I (g^{-1}(t)) e^{it\Delta }u_0\|_{X^{s,b}}$$
$$\lesssim |g(I)|^{1/2-b}\|u_0\|_{H^s},$$
where in the last line we applied \eqref{Xsb_est_1} in the standard case $g(t)=t$ to conclude \eqref{Xsb_est_1} in our case.

To prove \eqref{Xsb_est_2} we apply the previous strategy, that is, after a change of variables we get
$$	\left	\|  g'(t)\int_0^t g'(s)S(t,s)w(s)ds\right  \|_{X^{s,b}_g}
= 	\left	\|  \int_0^{t} e^{i(t-s')\Delta}w(g^{-1}(s'))ds' \right  \|_{X^{s,b}}$$
$$\lesssim \| w(g^{-1}(s'))\|_{X^{s,b-1}}= \| g'(t)w\|_{X^{s,b-1}_g},$$
where in the last line we applied \eqref{Xsb_est_2} in the standard case and the relation between $X^{s,b}$ and $X^{s,b}_g$ spaces.

Inequality \eqref{Xsb_est_3} follows again from the standard case, indeed,
	$$\| \chi_I (t) g'(t)  |u|^2u\|_{X^{s,b-1}_g}
	=	\| \chi_{g(I)} (t)   |u(g^{-1}(t))|^2u(g^{-1}(t))\|_{X^{s,b-1}}$$
	$$\lesssim \| \chi_{g(I)} (t) u(g^{-1}(t))\|_{X^{s,b'}}^2 \| \chi_{g(I)} (t) u(g^{-1}(t))\|_{X^{s,b}}$$
	$$=\| \chi_I (t) g'(t) u\|_{X^{s,b'}_g}^2 \| \chi_I (t) g'(t) u\|_{X^{s,b}_g},$$
	which proves \eqref{Xsb_est_3}.

More generally one has that the following trilinear estimate holds.

\begin{proposition}\label{prop.trilin.est}
	Let $s>0$, $b>1/2$ and $1/4<b'<b$. Let also $\chi_I$ be a smooth cutoff function as before. Then
\begin{equation} \label{trilinear-est}
\|g'(t) \,\,\chi_I v_1\,\, \chi_I v_2 \,\, \chi_I v_3\|_{X^{s,b-1}_g}
\lesssim \| \chi_I g'(t)v_1\|_{X^{s,b}_g} \| \chi_I g'(t)v_2\|_{X^{s,b'}_g} \| \chi_I g'(t)v_3\|_{X^{s,b'}_g} 
\end{equation}
\end{proposition}
\proof
By using the identity
$$ \|g'(t) \,\,\chi_I v_1\,\, \chi_I v_2 \,\, \chi_I v_3\|_{X^{s,b-1}_g}
=\|\,\,\chi_{g(I)} w_1\,\, \chi_{g(I)} w_2 \,\, \chi_{g(I)}w_3\|_{X^{s,b-1}},$$
with $w_i(t,x)=w_i(g^{-1}(t),x)$, $i=1,2,3$, the proof follows from the standard case $g(t)=t$ (see  \cite{Bourgain-Geom.Funct.Anal-1993}).
\endproof
Finally, \eqref{Xsb_est_4} can be proved again by means of the standard case, so we omit the details. This concludes the proof.
\endproof

We will now prove some results holding in the classical constant coefficients case from which the suitable formulation in our degenerate case will be derived. 

\begin{proposition}\label{prop.bilin1}
	Let $u\in X^{s,b}(\mathbb{R}\times \mathbb{T}^2)$ and let $f=f(t)$  be such that $f\in H^{1,b}(\mathbb{R})$. 
	Then, for all $s>0$ and for all $b\in (1/2,1)$, 
	\begin{equation}\label{bilin.1}
	\| u f\|_{X^{s,b}}\lesssim  \| f\|_{H^{1,b}}\|u\|_{X^{s,b}}.
	\end{equation}
\end{proposition}
\proof
Since
$$\| u f\|_{X^{s,b}}=\|\widehat{u\ast f}(\tau,k)\langle k\rangle^s \langle \tau-|k|^2\rangle^b\|_{L^2_t\ell^2_k}$$
we can prove the result by duality on $L^2_\tau\ell^2_k$.
We then take $v\in L^2_\tau\ell^2_k$ such that $\|v\|_{L^2_\tau\ell^2_k}=1$, and have
$$\sum_{k\in\mathbb{Z}^2}\int_{\mathbb{R}} |\widehat{u\ast f}|(\tau,k)\langle k\rangle^s \langle \tau-|k|^2\rangle^b |v|(\tau,k)d\tau$$
$$\lesssim\sum_{k\in\mathbb{Z}^2}\int_{\mathbb{R}^2}\langle\tau_1+\tau_2-|k|^2\rangle^{b}\langle k\rangle^s |\widehat{u}|(\tau_1,k)|\widehat{f}|(\tau_2) |v|(\tau_1+\tau_2,k)d\tau_1 d\tau_2$$
$$\lesssim  \sum_{k\in\mathbb{Z}^2}\int_{\mathbb{R}^2}\langle\tau_1-|k|^2\rangle^{b}\langle k\rangle^s |\widehat{u}|(\tau_1,k)|\widehat{f}|(\tau_2) |v|(\tau_1+\tau_2,k)d\tau_1 d\tau_2$$
$$+ \sum_{k\in\mathbb{Z}^2}\int_{\mathbb{R}^2}\langle\tau_2\rangle^{b}\langle k\rangle^s |\widehat{u}|(\tau_1,k)|\widehat{f}|(\tau_2) |v|(\tau_1+\tau_2,k)d\tau_1 d\tau_2$$
$$=I+II.$$
For the term $I$ we get
$$I=\sum_{k\in\mathbb{Z}^2}\int_{\mathbb{R}^2}\langle\tau_1-|k|^2\rangle^{b}\langle k\rangle^s |\widehat{u}|(\tau_1,k)|\widehat{f}|(\tau_2) |v|(\tau_1+\tau_2,k)d\tau_1 d\tau_2$$
$$\lesssim  \int_{\mathbb{R}_{\tau_2}} |\widehat{f}(\tau_2)|\left( \int_{\mathbb{R}_{\tau_1}}\| \langle\tau_1-|k|^2\rangle^{b}\langle k\rangle^s\widehat{u}(k) \|_{\ell^2_k}\|v(\tau_1+\tau_2,k)\|_{\ell^2_k}d\tau_1\right)d\tau_2$$
$$\lesssim \| \widehat{f}\|_{L^1}  \| u\|_{X^{s,b}}\|v\|_{L^2_\tau\ell^2_k}.$$

Similar computation on the term $II$ give
$$II=\sum_{k\in\mathbb{Z}^2}\int_{\mathbb{R}^2}\langle\tau_2\rangle^{b}\langle k\rangle^s \widehat{u}(\tau_1,k)\widehat{f}(\tau_2) v(\tau_1+\tau_2,k)d\tau_1 d\tau_2$$
$$\leq \|f\|_{H^{1,b}}\|u\|_{X^{s,b}}\|v\|_{L^2_\tau\ell^2_k}.$$
Finally, putting together the estimates for $I$ and $II$ the result follows.
\endproof

\begin{proposition}\label{prop-blin-space}
Let $\chi_I$ be a smooth cutoff function supported on $[-2\delta,2\delta]$ such that $\chi_I\equiv 1$ on $[-\delta,\delta]$, and $\beta \in H^{s+2b}(\mathbb{T}^2)$. Then
$$\|\chi_I \beta\|_{X^{s,b}}\lesssim \|\chi_I\|_{H_t^b}\|\beta\|_{H_x^{s+2b}}\lesssim \delta^{1/2-b}\|\beta\|_{H_x^{s+2b}}.$$
\end{proposition}
\proof
By definition of $X^{s,b}$ spaces and by the properties of the functions $\chi_I$ and $\beta$ we have
that
\begin{eqnarray}
	\|\chi_I\beta\|_{X^{s,b}}&=&\Big(\sum_{k\in\mathbb{Z}^2}\int_{\mathbb{R}}\langle k\rangle^{2s} \langle \tau -|k|^2\rangle^{2b} |\widehat{\chi_I}(\tau)|^2 |\widehat{\beta}(k)|^2 d\tau\Big)^{1/2} \nonumber\\
	&\lesssim& \Big(\sum_{k\in\mathbb{Z}^2}\int_{\mathbb{R}}\langle k\rangle^{2(s+2b)} \langle \tau\rangle^{2b}  |\widehat{\chi_I}(\tau)|^2 |\widehat{\beta}(k)|^2 d\tau\Big)^{1/2} \nonumber\\
	&\lesssim& \|\chi_I\|_{H_t^b}\|\beta\|_{H_x^{s+2b}}.\nonumber
\end{eqnarray}
We now write $\chi_I$ as $\chi_I(t)=\psi(\delta^{-1}t)$, with $\psi$ supported on $[-2,2]$ and $\psi\equiv 1$ on $[-1,1]$, and have,  since $\widehat{\chi_I}(\tau)=\delta \widehat{\psi}(\delta \tau)$, that
\begin{eqnarray}
\|\chi_I\|_{H^b_t}&=&\Big(\delta^2\int_{\mathbb{R}}\langle \tau \rangle^{2b}|\widehat{\psi}(\delta \tau)|^2d\tau \Big)^{1/2}\nonumber\\
&=& \Big(\delta\int_{\mathbb{R}}\langle \delta^{-1}\tau \rangle^{2b}|\widehat{\psi}( \tau)|^2d\tau\Big)^{1/2} \lesssim \delta^{1/2-b}\| \psi\|_{H^b_t}. \nonumber
\end{eqnarray}
Finally, combining the two estimates above we obtain the desired inequality.
\endproof

\begin{proposition}\label{quadrilin.est}
Let $h\in H^{1}(\mathbb{T}^2)$, $s>0$,  $b'<b$ and $b\in (1/2,1)$. Then, for $p_1>1/2$ and $s_1>1$,
$$\|\chi_I h\,  \chi_I u_1\, \chi_I u_2\,  \chi_I u_3\|_{X^{s,b-1}}\lesssim \|\chi_Ih\|_{H^{p_1}_tH^{s_1}_x}\|\chi_I u_1\|_{X^{s,b}}\|\chi_I u_2\|_{X^{s,b'}}\|\chi_I u_3\|_{X^{s,b'}}.$$
\end{proposition}
\proof
The proof follows by using the standard technique employed to prove $X^{s,b}$-multilinear estimates (see \cite{Bourgain-Geom.Funct.Anal-1993}). 

We proceed by duality (on $L^2_{\tau}\ell^2_k$) observing that 

$$\|\chi_I h\,  \chi_I u_1\, \chi_I u_2\,  \chi_I u_3\|_{X^{s,b-1}}=
\!\sup_{\substack{v\in L^2_\tau\ell^2_{,k}\\ \| v\|_{L^2_\tau\ell^2_{k}}=1}}
\Big|\int_{\mathbb{R}_{\tau_0}}\int_{\mathbb{R}_{\tau_1}}\int_{\mathbb{R}_{\tau_2}} \int_{\mathbb{R}_{\tau_3}}\widehat{\chi_I h }(\tau_0,k_0)\widehat{\chi_Iu}_1(\tau_1,k_1)\widehat{\chi_I u}_2(\tau_2,k_2)$$
$$ \times \,\widehat{\chi_I u}_3(\tau_3,k_3)
 \langle \tau_0+\tau_1+\tau_2+\tau_3-|k_0+k_1+k_2+k_3|^2\rangle^{b-1} \langle k_0+k_1+k_2+k_3 \rangle^s$$ 
 \begin{equation}\label{est.quad.1}
 \times \, \,v(\tau_0+\tau_1+\tau_2+\tau_3,k_0+k_1+k_2+k_3)d\tau_0\tau_1 d\tau_2 d\tau_3 \Big|.
 \end{equation}

We now take the Littlewood-Paley decomposition of $h,u_1, u_2$ and $u_3$, that is we write
$$h=\sum_{N_0\in 2^\mathbb{Z}} P_{N_0} h,\quad u_j=\sum_{N_j\in 2^\mathbb{Z}} P_{N_j} u_j,\quad j=1,2,3,$$
where $P_N$ is the Fourier multiplier such that $\widehat{P_Nu}(k)$ is supported in the region $\{k\in\mathbb{Z}^2; N/2\leq |k|\leq N\}$, and replace them in \eqref{est.quad.1} to get

$$\eqref{est.quad.1}= \sup_{\substack{v\in L^2_\tau\ell^2_{,k}\\ \| v\|_{L^2_\tau\ell^2_{k}}=1}} \sum_{\substack{N_0,N_1,N_2,N_3}}\sum_{\substack{ |k_0|\sim N_0\\|k_1|\sim N_1\\ |k_2|\sim N_2\\|k_3|\sim N_3 }}\int_{\mathbb{R}_{\tau_0}}\int_{\mathbb{R}_{\tau_1}}\int_{\mathbb{R}_{\tau_2}}\int_{\mathbb{R}_{\tau_3}} |\widehat{\chi_IP_{N_0}h}(\tau_0,k_0)||\widehat{\chi_IP_{N_1}u}_1(\tau_1,k_1)|
$$
$$ \langle \tau_0+\tau_1+\tau_2+\tau_3-|k_0+k_1+k_2+k_3|^2\rangle^{b-1} |\widehat{\chi_IP_{N_2}u}_2(\tau_2,k_2) |\widehat{\chi_IP_{N_1}u}_3(\tau_3,k_3)| $$
\begin{equation}
\label{est.quad.2}
\times \langle k_0+k_1+k_2+k_3\rangle^s |v(\tau_0+\tau_1+\tau_2+\tau_3,k_0+k_1+k_2+k_3)|  d\tau_0d\tau_1 d\tau_2 d\tau_3.
\end{equation}

We then write 
\begin{equation}
\label{eq.sum2}
\sum_{N_0,N_1,N_2,N_3}\sum_{\substack{|k_0|\sim N_0\\|k_1|\sim N_1\\ |k_2|\sim N_2\\|k_3|\sim N_3 }}
=\sum_{N_0,N_2,N_3}\sum_{N_1\sim  \max\{N_0,N_2,N_3\}}\sum_{\substack{ |k_0|\sim N_0\\|k_1|\sim N_1\\ |k_2|\sim N_2\\|k_3|\sim N_3  }}
+\sum_{N_0,N_2,N_3}\sum_{N_1\gg  \max\{N_0,N_2,N_3\}} \sum_{\alpha}\sum_{\substack{ k_1\in Q_\alpha\\ |k_0|\sim N_0\\|k_2|\sim N_2\\|k_3|\sim N_3 }},
\end{equation}
where $Q_\alpha$ is a cube of side $\sim  \max\{N_0,N_2,N_3\}$ and where $\alpha$ ranges over the number of cubes of side $\sim  \max\{N_0,N_2,N_3\}$ covering a cube of side $N_1$. 
From \eqref{est.quad.2} and \eqref{eq.sum2} we obtain
$$\|\chi_I h\,  \chi_I u_1\, \chi_I u_2\,  \chi_I u_3\|_{X^{s,b-1}}\lesssim
\sup_{\substack{v\in L^2_\tau\ell^2_{,k}\\ \| v\|_{L^2_\tau\ell^2_{k}}=1}}
\Bigg( \sum_{N_0,N_2,N_3,}\sum_{N_1\gg \max\{N_0,N_2,N_3\}} \sum_{\alpha} N_1^s $$
$$\int_{\mathbb{R}^4}|\widehat{ \chi_IP_{N_0}h}(\tau_0,k_0)|
|\widehat{\chi_I P_{Q_\alpha}u}_1(\tau_1,k_1)| |\widehat{\chi_IP_{N_2}u}_2(\tau_2,k_2)| |\widehat{\chi_IP_{N_3}u}_3(\tau_3,k_3)|$$
$$\langle\tau_0+\tau_1+\tau_2+\tau_3-|(k_0+k_1+k_2+k_3)|^2 \rangle^{b-1}|\widehat{P_{\tilde{Q}_\alpha}\check{v} }(\tau_0+\tau_1+\tau_2+\tau_3,k_0+k_1+k_2+k_3)| d\tau_0d\tau_1 d\tau_2 d\tau_3$$

$$+\!\!\!\!\sum_{N_0,N_2,N_3,}\sum_{N_1\sim \max\{N_0,N_2,N_3\} }\!\!\!\!\!\!\!\! N_1^s \int_{\mathbb{R}^4}|\widehat{ \chi_IP_{N_0}h}(\tau_0,k_0)|
|\widehat{\chi_I P_{Q_\alpha}u}_1(\tau_1,k_1)| |\widehat{\chi_IP_{N_2}u}_2(\tau_2,k_2)| |\widehat{\chi_IP_{N_3}u}_3(\tau_3,k_3)|$$
$$\langle\tau_0+\tau_1+\tau_2+\tau_3-|(k_0+k_1+k_2+k_3)|^2 \rangle^{b-1}|\widehat{P_{\tilde{Q}_\alpha}\check{v} }|(\tau_0+\tau_1+\tau_2+\tau_3,k_0+k_1+k_2+k_3)| d\tau_0d\tau_1 d\tau_2 d\tau_3\Bigg).$$

We now denote by $\check{v}$ the inverse Fourier transform of $v(\tau,k)$ and by $\tilde{Q}_\alpha$ a suitable cube of side $\sim N_2$.  With these notations in mind, by using  Plancherel's theorem and the property $|\widehat{u}\ast \widehat{\chi}_I |\leq |\widehat{u}|\ast \widehat{\chi}_I=(|\widehat{u}|^\vee\chi_I)^\wedge$ for $\chi_I$ such that $\widehat{\chi}_I\geq 0$ (for more details see \cite{Bourgain-VIM}), we have

$$\|\chi_I h\,  \chi_I u_1\, \chi_I u_2\,  \chi_I u_3\|_{X^{s,b-1}}\lesssim
\sup_{\substack{v\in L^2_\tau\ell^2_{,k}\\ \| v\|_{L^2_\tau\ell^2_{k}}=1}}
\Bigg( \sum_{N_0,N_2,N_3,}\sum_{N_1\gg \max\{N_0,N_2,N_3\}} \sum_{\alpha} N_1^s $$
$$\int_{\mathbb{R}\times\mathbb{T}^2}\chi_I(t)|\widehat{P_{N_0}h}|^\vee(x)
 |\widehat{P_{Q_\alpha}u}_1|^\vee(t,x)\times  |\widehat{P_{N_2}u}_2|^\vee(t,x) |\widehat{P_{N_3}u}_3|^\vee(t,x)$$
 $$(\langle\tau_0+\tau_1+\tau_2+\tau_3-|(k_0+k_1+k_2+k_3)|^2 \rangle^{b-1}|\widehat{P_{\tilde{Q}_\alpha}\check{v} }|)^{\vee}(t,x) \chi_I^3(t) \,dt dx$$
 $$+\sum_{N_0,N_2,N_3,}\sum_{N_1\sim \max\{N_0,N_2,N_3\} } N_1^s \int_{\mathbb{R}\times\mathbb{T}^2} \chi_I(t)|\widehat{P_{N_0}h}|^\vee(x)| |\widehat{P_{N_1}u}_1|^\vee(t,x) $$
 $$\times |\widehat{P_{N_2}u}_2|^\vee(t,x) |\widehat{P_{N_3}u}_3|^\vee(t,x)
  (\langle\tau_0+\tau_1+\tau_2+\tau_3-|(k_0+k_1+k_2+k_3)|^2 \rangle^{b-1}|\widehat{P_{\tilde{Q}_\alpha}\check{v} }|)^{\vee}\chi_I^3(t) \,dt dx\Bigg)$$
 $$\cong \sup_{\substack{v\in L^2_\tau\ell^2_{,k}\\ \|v\|_{L^2_\tau\ell^2_{k}}=1}}(I+II).$$

Now the delicate part to analyze is the one where the condition $N_1\gg \max\{N_0,N_3,N_3\}$ appears, namely, the first term $I$ on the right hand side of the inequality above.
For such a term, on denoting by $M:=\max\{N_0,N_3,N_3\}$, we have
$$I\lesssim \sum_{N_0,N_2,N_3,}\sum_{N_1\gg M} \sum_{\alpha} N_1^s
\|\chi_I(t)|\widehat{P_{N_0}h}|^\vee(x)\|_{L^\infty_{t,x}}
\|\chi_I|\widehat{P_{Q_\alpha}u}_1|^\vee\|_{L^4_{t,x}} \| \chi_I |\widehat{P_{N_2}u}_2|^\vee\|_{L^4_{t,x}}$$
$$
 \| \chi_I |\widehat{P_{N_3}u}_3|^\vee\|_{L^4_{t,x}}\|  (\langle\tau_0+\tau_1+\tau_2+\tau_3-|(k_0+k_1+k_2+k_3)|^2 \rangle^{b-1}|\widehat{P_{\tilde{Q}_\alpha}\check{v} }|)^{\vee}\|_{L^4_{t,x}}
$$

$$\lesssim \sum_{N_0,N_2,N_3,}\sum_{N_1\gg M} \sum_{\alpha} N_1^s
\|\chi_I(t)|\widehat{P_{N_0}h}|^\vee(x)\|_{H^{p_1}_tH^{s_1}_{x}}
\|\chi_I|\widehat{P_{Q_\alpha}u}_1|^\vee\|_{L^4_{t,x}} \| \chi_I |\widehat{P_{N_2}u}_2|^\vee\|_{L^4_{t,x}}$$
$$
\| \chi_I |\widehat{P_{N_3}u}_3|^\vee\|_{L^4_{t,x}}\|  (\langle\tau_0+\tau_1+\tau_2+\tau_3-|(k_0+k_1+k_2+k_3)|^2 \rangle^{b-1}|\widehat{P_{\tilde{Q}_\alpha}\check{v} }|)^{\vee}\|_{L^4_{t,x}}
$$
$$\lesssim \|\chi_I h\|_{H^{p_1}_tH^{s_1}_x}\|\chi_I u_1\|_{X^{s,b}}
\|\chi_I u_2\|_{X^{s,b'}}\|\chi_I u_3\|_{X^{s,b'}}\|\check{v}\|_{L^2_{t,x}},
$$
where in the first line we used Sobolev embeddings and H\"{o}lder's inequality, while in the last one we applied 
$$\|\chi_I P_N f\|_{L^4_{t,x}}\leq N^{\varepsilon}\|\chi_If\|_{X^{0,b_1}},\quad b_1\geq 1/2(1-\min\{\varepsilon,1/2\}),$$ 
 and choose $b'> 1/2(1-\varepsilon)$, with $\varepsilon$ satisfying $3\varepsilon<s$.  Finally, taking the supremum over the functions $v\in L^2_t\ell^2_k$ such that $\|v\|_{L^2\ell^2_k}=1$ we get the desired estimate for the term $I$. Finally, by similar arguments, the required estimate can be proved for the term $II$ as well, which concludes the proof.
\endproof

As a corollary of Proposition \ref{quadrilin.est} we get the inequality below.
\begin{corollary}\label{prop.bilin.2}
	Let $s>0$ and $b\in (1/2,1)$ Then, for all $s>0$, $b'<b$ and $b\in (1/2,1)$, we have
	$$\| \chi_I u_1 \chi_I u_2\|_{X^{s,b-1}}\lesssim \|\chi_Iu_1\|_{X^{s,b}}\|\chi_Iu_2\|_{X^{s,b'}}.$$
\end{corollary}
\proof
The proof follows form Proposition \ref{quadrilin.est} by taking $h=u_3=1$.
\endproof

From Proposition \ref{prop.trilin.est}, \ref{prop.bilin1}, \ref{prop-blin-space}, \ref{quadrilin.est}  and  Corollary \ref{prop.bilin.2}  we derive the suitable (weighted) formulation of the multilinear estimates holding in the time-degenerate case.

\begin{proposition}\label{bilin.g}
Let $s>0$, $b\in (1/2,1)$, $b'<b$, and $H^{p,b}_g(\mathbb{R})$ as in Definition \ref{def.Hspaces}. Then, for $h\in H^1(\mathbb{T}^2)$ and $\beta \in H^{s+2b}(\mathbb{T}^2)$, we have
\begin{equation}
\label{eq.bilin.g}
\|g'(t) f(t)u\|_{X^{s,b}_g}\lesssim \|g'f\|_{H^{1,b}_g}\|g'u\|_{X^{s,b}_g},
\end{equation}
\begin{equation}
\| g'(t)\chi_I \beta\|_{X^{s,b}_g}\lesssim \|g' \chi_I\|_{H^{2,b}_g}\|\beta\|_{H_x^{s+2b}},
\end{equation}
\begin{equation}
\label{eq.bilin.g2}
\|g'(t)\, \chi_I u_1 \, \chi_I u_2\|_{X^{s,b-1}_g}\lesssim \|g'(t) \chi_I u_1\|_{X^{s,b}_g}\|g'(t)\chi_Iu_2\|_{X^{s,b'}_g},
\end{equation}
and, for $p_1>1/2, s_1>1$,
$$\|g'(t)\, \chi_I h\, \chi_I u_1 \, \chi_I u_2 \, \chi_I u_3\|_{X^{s,b-1}_g}\lesssim \|g'(t)\chi_I h\|_{H^{2,p_1}_gH^{s_1}_x} \|g'(t) \chi_I u_1\|_{X^{s,b}_g}$$
\begin{equation}\label{eq.bilin.g3}
\times \|g'(t)\chi_Iu_2\|_{X^{s,b}_g}\|g'(t)\chi_Iu_2\|_{X^{s,b'}_g}.
\end{equation}
\end{proposition}
\proof
The proof follows immediately from Proposition \ref{prop.trilin.est}, \ref{prop.bilin1}, \ref{prop-blin-space}, \ref{quadrilin.est}  and  Corollary \ref{prop.bilin.2} by using the identities
$\|g'(t) v\|_{X^{s,b}_g}=\| v(g^{-1}(\cdot),\cdot)\|_{X^{s,b}}$ and $\widehat{v\circ g^{-1}}(\tau)= \widetilde{g' v}(\tau)$.
\endproof
Let us remark that we will not be using all the estimates in Proposition \ref{bilin.g} in the time-degenerate setting. However, since the estimates in Proposition \ref{prop-blin-space}, \ref{quadrilin.est} and in Corollary \ref{prop.bilin.2} involving standard $X^{s,b}$-spaces will be used in the next session, we translated them in the time-degenerate setting for completeness.

\subsection{Local well-posedness}
With the previous estimates at our disposal we can now focus on the proof of the local well-posedness of the IVP 
\begin{equation*}\label{IVPtime}
	\left\{
	\begin{array}{l}
	i\partial_t u + g'(t)\Delta_x u= g'(t)|u|^2u.\\
	u(0,x)=u_0(x)
	\end{array}
	\right.
\end{equation*}
The proof is centered on the contraction mapping theorem and by now the argument  is standard \cite{Bourgain-Geom.Funct.Anal-1993, Kenig-Ponce-Vega_JAMS_19996}. Still we prefer to report the details below since we will be using the modified spaces $\tilde{X}^{s,b}_g$. 

%

\proof[Proof of Theorem \ref{thm.time-lwp}]
We start by defining the metric space $X$ as
$$X:=\{u\in I\times \mathbb{T}^2\rightarrow \mathbb{C}; \| u\|_{\tilde{X}^{s,b}_g}<\infty\},$$
where  $\tilde{X}^{s,b}_g$ is as in Definition \ref{def.Xsbg}.


Now, given $u_0\in H^{s}(\mathbb{T}^2)$, and $I=[-2\delta, 2\delta]$, with $0<\delta<1$ to be determined later, we define the operator
$$\Phi_{u_0}( u):= \chi_I S(t)u_0+  \chi_I(t) \int_0^t S(t,t') g'(t') |\chi_I(t')u(t')|^2 \chi_I(t')u(t') dt'$$
where $\chi_I$ is a smooth cutoff function such that $\chi_I\equiv 1$ on $[-\delta,\delta]$, and prove that $\Phi_{u_0}$ is a contraction on $B_R:=\{u; \| u\|_{\tilde{X}^{s,b}_g}\leq R\}\subset \tilde{X}^{s,b}_g$, with $R= 2C'\delta^{1/2-b}\|u_0\|_{H^s(\mathbb{T}^2)}$.

Note that, by using Proposition \ref{prop.xspEst} and the fact that $|b(I)|\lesssim \delta$ (which is due to the properties of $b$), we get, for $u\in B_R$,
$$\| \Phi_{u_0}(u)\|_{\tilde{X}^{s,b}_g}=\| g'\Phi_{u_0}(u)\|_{X^{s,b}_g}$$
$$=\|  g'(t)\chi_I S(t)u_0\|_{X^{s,b}_g}+ \| g'(t) \chi_I(t) \int_0^t S(t,t') g'(t') |\chi_I(t')u(t')|^2 \chi_I(t')u(t') dt'\|_{X^{s,b}_g}$$
$$\lesssim  |b(I)|^{1/2-b}\| u_0\|_{H^s(\mathbb{T}^2)}+ \|\chi_I g'(t)u\|^2_{X^{s,b'}_g}  \|\chi_I g'(t)u\|_{X^{s,b}_g}$$
$$\leq C \delta^{1/2-b}\| u_0\|_{H^s(\mathbb{T}^2)}+ C^3 \delta^{1/2-b}\delta^{(b-b')/4} \|\chi_I g'(t)u\|_{X^{s,b}_g}^3$$
$$\leq \frac R 2+  \delta^{1/2-b}\delta^{(b-b')/4}R^3.$$
Since $\frac 1 4 < b'< b$, $b\in (1/2,1)$, then, taking $b=1/2 +\varepsilon$, with $\varepsilon\in (0,1/16)$ (so $b\in (\frac 1 2, \frac {9}{16})$), and $b'= 1/2-4\varepsilon>1/4$, we have that $1/2-b+(b-b')/4=\varepsilon/4>0$. Therefore, by choosing $\delta$ sufficiently small and such that $\delta^{1/2-b}\delta^{(b-b')/4}R^2\leq 1/2$, from the previous inequality we obtain that $\Phi_{u_0}$ maps the ball $B_R$ into itself.


To conclude that $\Phi_{u_0}$ is a contraction we write
$$\chi_I u |\chi_I u|^2-\chi_I v |\chi_I v|^2=\chi_I (u-v) |\chi_I u|^2+\chi_Iv \chi_I\overline{u}\chi_I(u-v)
+|\chi_I v|^2\chi_I(\overline{u}-\overline{v})$$
and apply \eqref{trilinear-est} to get 
$$\| \Phi_{u_0}(u)-\Phi_{u_0}( v)\|_{\tilde{X}^{s,b}_g}
\lesssim (\| \chi_I g'(t) u\|_{X^{s,b'}_g}^2+ \| \chi_I g'(t) v\|^2_{X^{s,b'}_g} )  \| \chi_I g'(t) (u-v)\|_{X^{s,b}_g}$$
$$\lesssim R^3  \delta^{1/2-b+(b-b')/4}\| u-v\|_{\tilde{X}^{s,b}_g}<\| u-v\|_{\tilde{X}^{s,b}_g},$$
where the last inequality follows by choosing, eventually, $\delta$ smaller than before. 
This, finally, gives the result for $t\in  [-\delta,\delta]$, since, recall,  $\chi_I\equiv 1$ on  $[-\delta,\delta]$, and
$$u(t,x)=S(t)u_0+  \int_0^t S(t,t') g'(t') |u(t')|^2 (t')u(t') dt',\quad \text{for all}\,\,  t\in[-\delta,\delta].$$
\endproof

By using the previous approach we can also prove the local well-posedness of the IVP
\begin{equation*}\label{IVPtime2}
\left\{
\begin{array}{l}
i\partial_t u + g'(t)\Delta_x u= f(t)|u|^2u,\\
u(0,x)=u_0(x),
\end{array}
\right.
\end{equation*}
on $I\times \mathbb{T}^2$, where $I$ is a suitable finite interval of time and $f$ is smooth enough, namely $f\in H^{1,b}_g$ for $b\in (1/2,1)$.  Note that assuming $f\in H^{1,b}_g$ in particular implies  that $f$ has to be zero at time zero of the same order of $g'(t)$.

\proof[Proof of Theorem \ref{thm.time-lwp2}]
We start by defining the metric space $X$ as
$$X:=\{u\in I\times \mathbb{T}^2\rightarrow \mathbb{C}; \| u\|_{\tilde{X}^{s,b}_g}<\infty\},$$
where 
$$\tilde{X}^{s,b}_g:=\{u\in X^{s,b}_g; g' u\in X^{s,b}_g\},$$
and
$$\| u\|_{\tilde{X}^{s,b}_g}=\| g' u\|_{X^{s,b}_g}.$$

Now, given $u_0\in H^{s}(\mathbb{T}^2)$ such that $\|u_0\|_{H^s}=r$, and $I=[-2\delta, 2\delta]$, with $0<\delta<1$, we define the operator
$$\Phi_{u_0}(u):= \chi_I S(t)u_0+  \chi_I(t) \int_0^t S(t,t') f(t') |\chi_I(t')u(t')|^2 \chi_I(t')u(t') dt'$$
where $\chi_I$ is a smooth cutoff function such that $\chi\equiv 1$ on $[-\delta,\delta]$, and prove that $\Phi_{u_0}$ is a contraction on $B_R:=\{u; \| u\|_{\tilde{X}^{s,b}_g}\leq R\}\subset \tilde{X}^{s,b}_g$, with $R= 2C'\delta^{1/2-b}\|u_0\|_{H^s(\mathbb{T}^2)}$.

Note that, by using Proposition \ref{prop.xspEst}, Proposition \ref{bilin.g}, and the fact that $|b(I)|\lesssim \delta$ (which is due to the properties of $b$), we get, for $u\in B_R$,
$$\| \Phi_{u_0}(u)\|_{\tilde{X}^{s,b}_g}=\| g'\Phi_{u_0}(u)\|_{X^{s,b}_g}$$
$$=\|  g'(t)\chi_I S(t)u_0\|_{X^{s,b}_g}+ \| g'(t)  \int_0^t S(t,t') g'(t')  \frac{f(t')}{g'(t')} |\chi_I(t')u(t')|^2 \chi_I(t')u(t') dt'\|_{X^{s,b}_g}$$
$$\lesssim  |b(I)|^{1/2-b}\| u_0\|_{H^s(\mathbb{T}^2)}+ \|\chi_I g'(t)u\|^2_{X^{s,b'}_g}  \|g'(t')\chi_I u \,  \frac{f(t')}{g'(t')}\|_{X^{s,b}_g}$$
$$\lesssim  |b(I)|^{1/2-b}\| u_0\|_{H^s(\mathbb{T}^2)}+ \|\chi_I g'(t)u\|^2_{X^{s,b'}_g}  \|\chi_I g'(t) u \|_{X^{s,b}_g}  \|g'(t)\frac{f(t)}{g'(t)}\|_{H^{1,b}_g}$$

$$\leq C \delta^{1/2-b}\| u_0\|_{H^s(\mathbb{T}^2)}+ C^3 \delta^{1/2-b}\delta^{(b-b')/4}\| g'u\|_{X^{s,b}_g}^3$$
$$\le \frac R 2+ \delta^{1/2-b}\delta^{(b-b')/4} R^3\leq R$$
for  $b, b'$ as before and $\delta$ sufficiently small.
To prove that $\Phi_{u_0}$ is a contraction we follow the same steps as before. This completes the proof.
%
\endproof
\medskip

\textbf{The one dimensional case.} The previous techniques apply with few suitable modifications in the case when the problem is set on $\mathbb{R}\times \mathbb{T}$. Here to prove that the  quintic semilinear problem \eqref{IVPtime2Intro-Q} is locally well-posed with solution (multiplied by a suitable cutoff function in time) belonging to the $X^{s,b}_g(\mathbb{R}\times \mathbb{T})$ spaces, one translates, as in the previous case, the well-known standard Strichartz and multilinear estimates into suitable weighted Strichartz and multilinear estimates. We omit the details which are left to the interested reader.

\proof[Proof of Theorem \ref{thm.time-lwp-Q}]
The proof is straightforward after applying the technique used in the two-dimensional case based on the use of $X^{s,b}_g$-spaces.
\endproof

\section{A space-variable coefficients case}\label{Section.space.variable}
This section is devoted to the study of the of the IVP
\begin{equation}\label{IVP3}
\left\{	\begin{array}{l}
i\partial_t u+a_1(x_1)\partial_{x_1}^2 u+a_2(x_2)\partial_{x_2}^2 u=  u|u|^2,\\
u(0,x)=u_0(x),
\end{array}\right.
\end{equation}
with $a_1,a_2\in C^\infty (\mathbb{T})$ real valued and strictly positive. As before, we treat in detail the two dimensional case above, since, with few modifications, the result in Theorem \ref{thm.lwp.space-IntroQ} about the quintic one-dimensional case can be proved.

Our strategy to treat this problem consists in combining  a change of variables  with a gauge transformation  to reduce the problem to one where the linear part of the equation has constant coefficients.

We recall that a gauge transformation  is a multiplication operator of the form
$$Tf(t,x):=e^{\Phi(t,x)}f(t,x)$$
where the function $\Phi$ (in this periodic setting) is periodic in $x$.

We start by applying in \eqref{IVP3} the change of variables
$$(x_1,x_2)=(\alpha_1(y_1),\alpha_2(y_2)):=\alpha(y),$$
so that, on denoting by $v(t,y):=u(t,\alpha(y))$ and by $v_0(y):=u_0(\alpha(y))$, and assuming that $u$ solves \eqref{IVP3}, we have
\begin{eqnarray}
 i\partial_t v(t,y)+\Delta_y v(t,y)&=& i\partial_tu(t,\alpha(y))+(\partial^2_{x_1}u)(t,\alpha(y))(\partial_{y_1}\alpha_1(y_1))^2 +(\partial^2_{x_2}u)(t,\alpha(y))(\partial_{y_2}\alpha_2(y_2))^2\nonumber\\ &+&(\partial_{x_1}u)(t,\alpha(y))\partial^2_{y_1}\alpha_1(y_1)
 +(\partial_{x_2}u)(t,\alpha(y))\partial^2_{y_2}\alpha_2(y_2).
\nonumber
\end{eqnarray}
Then, by choosing $\alpha(y)=(\alpha_1(y_1),\alpha_2(y_2))$ such that
$\partial_{y_1}\alpha_1(y_1)=\sqrt{a_1(\alpha_1(y_1))}$ and $\partial_{y_2}\alpha_2(y_2)=\sqrt{a_2(\alpha_2(y_2))}$, and using the fact that $(\partial_{x_j}u)(t,\alpha(y))=(\partial_{y_j}v(t,y))\partial_{y_j}\alpha_j(y_j)$ for all $j=1,2$, we get that $v$ solves
\begin{equation}\label{newEq}
\left\{\begin{array}{l}
 i\partial_t v(t,y)+\Delta_yv(t,y)-(\partial_{y_1}v(t,y))\frac{\partial^2_{y_1}\alpha_1(y_1)}{\partial_{y_1}\alpha_1(y_1)} 
-(\partial_{y_2}v(t,y))\frac{\partial^2_{y_2}\alpha_2(y_2)}{\partial_{y_2}\alpha_2(y_2)} =v|v|^2\\
v(0,y)=v_0(y)
\end{array}\right.
\end{equation}


We now apply the gauge transform 
$$Tf(t,y):=e^{\Phi(y)}v(t,y)=\exp\Big\{-\frac 1 2\int_{0}^{y_1}\frac{\alpha_1''(s_1)}{\alpha_1'(s_1)}ds_1-\frac 1 2 \int_{0}^{y_2}\frac{\alpha_2''(s_2)}{\alpha_2'(s_2)}ds_2 \Big\}f(t,y),$$
(where, note, $\Phi$ is periodic since $\partial_{y_j}\alpha_j(y_j)=\sqrt{a_j(\alpha_j(y_j))}$) on the left and right hand side of \eqref{newEq} and have, on denoting by $w(t,y):=e^{\Phi(y)} v(t,y)$,
$$
i\partial_t w+\Delta_y w-
\left( \partial_{y_1}^2\Phi+ \partial_{y_2}^2\Phi+(\partial_{y_1}\Phi)^2+(\partial_{y_2}\Phi)^2\right)w=e^{-2\Phi}w|w|^2.$$

Finally, we reduced the study of \eqref{IVP3} to the study of
\begin{equation}\label{IVP4}
\left\{\begin{array}{l}
i\partial_t w+\Delta_y w=e^{-2\Phi}w|w|^2-\beta w\\
w(0,y)=w_0(y)
\end{array}\right.
\end{equation}
with $\beta=\beta(y)=\partial_{y_1}^2\Phi+ \partial_{y_2}^2\Phi+(\partial_{y_1}\Phi)^2+(\partial_{y_2}\Phi)^2$, and 
$w_0(y)= e^{\Phi(y)} u_0(\alpha(y))$.

\begin{theorem}\label{thm.lwp.space}
Let $s>0$ and $\beta\in H^{s+2b'}(\mathbb{T}^2)$ with $b'>1/2-\varepsilon$ for some $\varepsilon>0$. Then, given $\Phi$ such that $e^{-2\Phi}\in H^1(\mathbb{T}^2)$, for every $w_0\in H^s(\mathbb{T}^2)$ there exists a unique solution of the IVP \eqref{IVP4} in the time interval $[-T,T]$ for a suitable $T=T(\|w_0\|_{H^s(\mathbb{T}^2)})$. Moreover, the solution $w$ satisfies
$$w\in C([-T,T];H^s(\mathbb{T}^2))$$
and
$$\chi_Iw\in X^{s,b}(\mathbb{R}\times\mathbb{T}^2),$$
with $b\in (1/2,1)$ such that $b>b'$, $I$ a closed neighborhood of $[-T,T]$, and $\chi_I$ a smooth cutoff fuction such that $\chi_I\equiv 1$ on $[-T,T]$.

\end{theorem}
\proof
In this proof we consider the metric space $X=X^{s,b}( \mathbb{R}\times\mathbb{T}^2)$ and take $w_0\in H^s(\mathbb{T}^2)$. We then consider $I=[-2\delta,2\delta]$, with $0<\delta<1$ to be determined later, and define the operator
$$\Psi_{w_0}(w):=\chi_I S(t)w_0+\chi_I(t)\int_0^tS(t,t')\Big( \chi_I^2e^{-2\Phi}w(t')|\chi_Iw(t')|^2-\chi_I \beta(y) \chi_I w(t') \Big)dt',$$
where $\chi_I$ is a smooth cutoff function such that $\chi_I\equiv 1$ on $[-\delta,\delta]$.
Now, given $B_R:=\{w;\|w\|_{X^{s,b}}\leq R\}$, with $R=2C\delta^{1/2-b}\|w_0\|_{H^s(\mathbb{T}^2)}$ and $C$ suitable (see the rest of the proof), we prove that $\Psi_{w_0}$ is a contraction. By using Proposition \ref{prop.bilin1}, \ref{quadrilin.est} and Corollary \ref{prop.bilin.2}, we get, for $w\in B_R$, $p_1>1/2, s_1>1$ (here we take $p_1=1/2+\varepsilon'$, with $0<\varepsilon'<1$ small enough),
\begin{eqnarray}
\|\Psi_{w_0}(w)\|_{X^{s,b}}&\leq& C_1\delta^{1/2-b}\|w_0\|_{H^s(\mathbb{T}^2)}+C_2\|\chi_Ie^{-2\Phi}\|_{H^{p_1}_tH^{s_1}_x}\delta^{1/2-b+(b-b')/4}\|\chi_I w\|^3_{X^{s,b}}\nonumber\\
&+&  \|\chi_I \beta\|_{X^{s,b'}}\|\chi_Iw\|_{X^{s,b}} ,\nonumber\\
&\leq&C_1R+C_2'R+C'_3\delta^{1/2-b+(b-b')/4}R^3\leq R\label{contraction.est}
\end{eqnarray}
for $b,b'$ as above (see, for instance, the proof of Theorem \ref{thm.time-lwp}),  $\delta$ small enough and $C$ suitable, for instance $C=3\max\{C_1,C_2',C_3'\}$.  
Note that in \eqref{contraction.est} we used the estimate  in Proposition \ref{prop-blin-space}, that is, $ \|\chi_I \beta\|_{X^{s,b'}}\lesssim \|\chi_I\|_{H^{b'}_t}\|\beta\|_{H^{s+2b'}_y(\mathbb{T}^2)}\lesssim \delta^{1/2-b'} C\lesssim C$ for $b'=1/2-\varepsilon$, since $1/2-\varepsilon <b'<b, b\in (1/2,1)$.
By similar arguments, eventually by choosing a different suitable constant $C$, we can conclude that $\Phi_{w_0}$ is a contraction on $B_R$. This, finally, concludes the proof (with $T=\delta$).
\endproof

\proof[Proof of Theorem \ref{thm.lwp.space-Intro}]
In order to prove Theorem \ref{thm.lwp.space-Intro} it is enough to apply Theorem \ref{thm.lwp.space} to the IVP \eqref{IVP4}.
To do that we need that $\beta\in H^{s+b'}(\mathbb{T}^2)$ and that $e^{-2\Phi}\in H^1(\mathbb{T}^2)$. Since we have  that $a_1, \, a_2 \in C^\infty(\mathbb{T}^2)$, $a_1, \, a_2$ are  strictly positive,  $\alpha_1'(y_1)=\sqrt{a_1(\alpha_1(y_1))}>0$, and  $\alpha_2'(y_2)=\sqrt{a_2(\alpha_2(y_2))}>0$, then we get that $\Phi\in C^\infty(\mathbb{T}^2)$, and, consequently, that $\beta$ and $e^{-2\Phi}$ have the required properties.
Therefore, by Theorem \ref{thm.lwp.space}, the IVP \eqref{IVP4} is locally-well posed in $H^s$ for all $s>0$, and, going back to the IVP \eqref{newEq} with the inverse gauge transform, we can conclude that the latter  is also locally well posed in $H^s$, for all $s>0$, and that the solution $v$ satisfies $\chi_I v\in X^{s,b}_\Phi(\mathbb{R}\times \mathbb{T}^2)$, where, recall,
$$X^{s,b}_\Phi(\mathbb{R}\times \mathbb{T}^2):=\{f:\mathbb{R}\times \mathbb{T}^2\rightarrow \mathbb{C}; e^{\Phi}f\in X^{s,b}(\mathbb{R}\times \mathbb{T}^2)\}.$$
Finally, the solution $u$ of \eqref{IVP3},  will belong to the space $X^{s,b}_{\Phi,\tilde{\alpha}}:=\{f:\mathbb{R}\times \mathbb{T}^2\rightarrow \mathbb{C}; (e^{\Phi}\,f)\circ \tilde{\alpha}\in X^{s,b}(\mathbb{R}\times \mathbb{T}^2)\}$, where $\tilde{\alpha}: \mathbb{R}\times \mathbb{T}^2\rightarrow \mathbb{R}\times \mathbb{T}^2$ is the diffeomorphism given by $\tilde{\alpha}(t,y):=(t, \alpha(y))$.
This concludes the proof.
\endproof
\begin{remark}
Going back to the proof that we just concluded it is easy to see that   one does not need to require $C^\infty (\mathbb{T})$ regularity for the  coefficients $a_1,\, a_2$, in fact $H^2$ is enough. Here we did not want to concentrate on the regularity of the coefficients, but rather on the  required regularity of the initial data, and show that  it is the same as the one required for a completely flat torus, that is $H^s, \, s>0$.
\end{remark}
\medskip

\begin{remark}
By combining the arguments applied in this and in the previous sections one can also consider a problem of the form
\begin{equation}\label{IVP.mixed}
\left\{\begin{array}{l}
i\partial_tu+ g'(t)\sum_{j=1}^2a_j(x_j)\partial_{x_j} u= f(t)u|u|^2,\\
u(0,x)=u_0(x),
\end{array}\right.
\end{equation}
where $g,f$ and $a_j$, $j=1,2$, satisfy the same properties as above.

In this case, performing a change of variables in space first, and, afterwards, applying a gauge transform with a suitable $\Phi=\Phi(y)$,  one finally reduces the problem \eqref{IVP.mixed} to one of the form \eqref{IVPtime2} (eventually with a nonlinearity of the form $f(t)e^{-2\Phi}w|w|^2-g'(t)\beta(x)w$) to which Theorem \ref{thm.time-lwp2} applies with few modifications. Finally one goes back to the solution of the original initial value problem by using the inverse gauge transform and changing variables again. The solution obtained this way will belong to a suitable weighted space $\tilde{X}^{s,b}_{g,\Phi, \tilde{ \alpha}}(\mathbb{R}\times \mathbb{T}^2):=\{f:\mathbb{R}\times \mathbb{T}^2\rightarrow \mathbb{C}; (e^{\Phi}f)\circ \tilde{\alpha}\in \tilde{X}^{s,b}_g(\mathbb{R}\times \mathbb{T}^2)\}$ with $\tilde{ \alpha}$ as before.
\end{remark}
\medskip

\noindent\textbf{The one dimensional case.} Is it easy to see that in the one dimensional case, that is when the problem is considered on $\mathbb{R}\times\mathbb{T}$, the presence of a single space variable allows us to apply the same scheme applied in the two-dimensional case. 

\proof[Proof of Theorem \ref{thm.lwp.space-IntroQ}] The proof is straightforward and follows from the two-dimensional strategy. Therefore, after a change of variable in space followed by the application of a gauge transform with a space-dependent function $\Phi$, the problem is translated in a constant coefficients one.  The solution of the resulting constant coefficients nonlinear problem will, finally,  provide the solution of the original one via the use of the inverse gauge transform and a change of variables.

\begin{remark}
We remark again that we expect that all the results of this paper can be generalized to the higher dimensional setting with few modification. We also expect the results in the time-degenerate case to be true on general closed manifolds combining our strategy with the known results mentioned in the introduction.
\end{remark}

%

	
%
	
	\nocite{*}
	\bibliographystyle{unsrt}
	\bibliography{References}	
	
\end{document}